\magnification=\magstep1

\input amssym.def

\input amssym.tex

%

%

%\font\sevensl=cmsl7

%\font\fivesl=cmsl5

%\scriptfont5=\sevensl

%\scriptscriptfont5=\fivesl

%

\font\teneusm=eusm10

\font\seveneusm=eusm7

\font\fiveeusm=eusm5

\newfam\eusmfam

\textfont\eusmfam=\teneusm

\scriptfont\eusmfam=\seveneusm

\scriptscriptfont\eusmfam=\fiveeusm

\font\tenmib=cmmib10

\font\sevenmib=cmmib7

\font\fivemib=cmmib5

\newfam\mibfam

\textfont\mibfam=\tenmib

\scriptfont\mibfam=\sevenmib

\scriptscriptfont\mibfam=\fivemib

\font\tenss=cmss10

\font\sevenss=cmss8 scaled 833

\font\fivess=cmr5

\newfam\ssfam

\textfont\ssfam=\tenss

\scriptfont\ssfam=\sevenss

\scriptscriptfont\ssfam=\fivess

\thinmuskip = 2mu

\medmuskip = 2.5mu plus 1.5mu minus 2.1mu  %minus 2.5mu

\thickmuskip = 4mu plus 6mu

\def\loosegraf#1\par{{%

\baselineskip=13.4pt plus 1pt \lineskiplimit=1pt \lineskip=1.3 pt

#1\par}}

\def\cc{{\Bbb C}}

\def\eps{\varepsilon}

%\endinput

\magnification=\magstep1

\def\longpoints{\leaders\hbox to 0.5em{\hss.\hss}\hfill \hskip0pt}

\bigskip

\bigskip

\centerline{\bf The Hartogs-type extension theorem for  meromorphic
mappings into}

\smallskip

\centerline{\bf $q$-complete complex spaces }

\bigskip

\centerline{Sergei Ivashkovich,  Alessandro Silva}

\bigskip\noindent\bf 0. Introduction.

\rm The aim of this note is to prove a  result on extension of
meromorphic mappings, which can be considered as a direct generalisation
of the Hartogs extension theorem for holomorphic functions.

\medskip\rm

Let $\Delta_r^n$ be the polydisc of radius $r$ in $\bf C^n$, 
and set $\Delta^n:=\Delta_1^n$. Let us define the "$q$-concave"
Hartogs figure $H_n^q(r)$ as the following open set in $\bf C^{n+q}$:

$$
H_n^q(r):=\Delta^n\times (\Delta^q\setminus \bar\Delta_{1-r}^q)\cup
\Delta_r^n\times \Delta^q.\eqno(1)
$$

Note that  $H_n^q(r)$ has $\Delta^{n+q}$ as its envelope of holomorphy.

\smallskip

Let $Y$ be a reduced complex space. Meromorphic mappings with values in
$Y$ are sad to satisfy a Hartogs-type extension Theorem if any such
$f:H_n^q(r)\longrightarrow Y$  extends to a meromorphic map
map $\hat f:\Delta^{n+q}\to Y$ from the unit polydisk
$\Delta^{n+q}$ into $Y$.  Sometimes we shall say more precisely that
meromorphic maps into such $Y$ possed a meromorphic extension property
in bidimension $(n,q)$.

Hartogs-type extension Theorems for
meromorphic mappings have been proved when $Y$ is compact K\"ahler, and
when $Y$ has some weaker metric properties by the first author in [I$_1$],
[I$_2$] and [I$_3$]. In this note we shall prove the following

\bigskip\noindent\bf

Theorem. \it Every meromorphic mapping $f:H_n^q(r)\longrightarrow Y$,
where $Y$ is a $q$-complete complex space, extends to a meromorphic mapping
 from $\Delta^{n+q}$ to $Y$.

\rm \bigskip

We recall that a \it strictly $q$-convex \rm function $\rho$   on the
complex space $Y$ with ${\sl dim} Y=N$ is a real valued $C^2$ function such
that the hermitian matrix consisting of the coefficients of the $(1,1)$ real
form $dd^c\rho $ has at least $N-q+1$ positive eigenvalues at all points
of $Y$. (Smooth objects on a  complex space $Y$ are by definition
the pull-backs of smooth objects in domains of $\bf C^M$ \rm under
appropriate local embeddings. The number $q$ is independent of
such embeddings).

The complex space $Y$ is  called \it $q$-complete \rm if there
exists a strictly $q$-convex exhaustion function $\rho :Y\longrightarrow
R^+$.

\smallskip

 We remark that in the case $q=1$, i.e. when $Y$ is Stein, the statement
 of the  Theorem, via proper embedding of $Y$ into $\bf C^M$ \rm, reduces
 to the extension of  holomorphic functions. This is given by the classical
 theorem of Hartogs, [H].

 More generally our Theorem provides Hartogs'
 type extension of meromorphic mappings into a complex subspaces of
 $\bf CP^{ N} \setminus\bf CP^{ N-q}$, see paragraph 3. Note that the Stein
 case inkludes here as $\bf CP^{ N} \setminus\bf CP^{ N-1}$.

 One more point, which we would like to mention in the Introduction is that
 our Theorem improves the following result due to K. Stein:

 \smallskip\it

 Let $D$ be a domain in $\cc^{q+2}, q\ge 1$ and $K\subset\subset D$ a compact
 subset in $D$ with connected complement. Let further $Y$ be a normal complex
 space of dimension $q$. Then every holomorphic  mapping
 $f:D\setminus K\to Y$ extends holomorphically  onto $D$.

 \smallskip\rm

 Note now that every noncompact irreducible complex space of dimension
 $q$ is
 $q$-complete, see [O]. So we have the folloving immediate corollary
 from our Theorem:

 \medskip\noindent\bf

 Corollary 1. \it Let $Y$ be an irreducible non compact    complex
 analytic space of dimension $q$. Every meromorphic mapping
 $f:H_n^q(r)\longrightarrow Y, n\ge 1$, extends to a meromorphic mapping from
 $\Delta^{n+q}$ to $Y$.\rm

 \smallskip

 If $Y$ is compact and $f:H_n^q(r)\longrightarrow Y$ is not
 surjective, then delete from $Y$ one point, which is not in the image
 of $f$, and call the resulting space
 $Y'$. Corollary 1 now applies to  $f:H^q_n(r)\to Y'$ and $f$  again extends
 onto the corresponding polydics.

 %Using this observations we shall prove in \S 4 the following improvement

 %of the theorem of Stein:

 \smallskip

 In particular,

 \it for any domain $D\subset \cc^{q+1}$ and any compact
 $K\subset\subset D$ with connected complement, any nonsurjective
 meromorphic mapping   $f:D\setminus K \to Y$ extends meromorphically onto
 $D$.

 \smallskip\rm In our problem section we shall discuss among other open
 questions also some ones arising from the attempts to remove the condition
 on $Y$ to be noncompact in the last statement.

 \bigskip\bf

 Table of content

 \medskip\noindent

 \line{0. Introduction. \longpoints pp. 1 - 2}

 \smallskip\noindent

 \line{1. Pleminaries. \longpoints pp. 2 - 5}

 \smallskip\noindent

 \line{2. Proof of the Theorem. \longpoints pp. 5 - 7}

 \smallskip\noindent

 \line{3. Consequences and open questions. \longpoints pp. 7 - 8}

 \smallskip\noindent

 \line{References}

 \bigskip\noindent\bf 1. Preliminaries.

 \smallskip \rm

 Let $X$ and $Y$ be reduced complex spaces with $X$  normal. A \it
 meromorphic mapping \rm $f:X\longrightarrow  Y$ is defined as an
 irreducible,  locally irreducible analytic subset
 $\Gamma_f\subset X\times Y$ (the graph of $f$), such that the restriction
 to $\Gamma_f$,  $\pi\mid_{\Gamma_f}:\Gamma_f\longrightarrow X$,  of the
 natural projection $\pi:X\times Y\longrightarrow X$   is proper,
 surjective and generically one to one, see [R]. The set
 $f[x]:=\{ y\in Y:(x,y)\in \Gamma_f\} $ is a compact subvariety in $Y$.
 The  set of points $x\in X$ such that ${\sl dim}f[x]\ge 1$ is analytic by
 the Remmert proper mapping theorem  and has codimension at least two, because 
 of the condition of irreducibility of $\Gamma_f$.
 This set is called the \it fundamental set \rm of $f$  or the \it set  of
 points of indeterminancy \rm of $f$ and will be denoted by $F$. If $X_1$
 is a normal subspace of $X$, $X_1\not\subset F$, we denote by
 $f\mid_{X_1}$ the meromorphic mapping with a graph equal to the
 (unique!)  irreducible component  of $\Gamma_f\cap (X_1\times Y)$, which
 projects \it onto \rm $X_1$.

 We shall list now some statements needed for the proof of our Theorem.
 First of all let us define the set
 
 $$
 E_n^q(r) = (\Delta^{n-1}\times \Delta_r\times \Delta^q )\cup (\Delta^{n-1}
 \times \Delta \times A^q(1-r,1) = \Delta^{n-1}\times H^q_1(r).\eqno(2)
 $$

 Here $A^q(1-r,1):=\{ z\in \cc^q: 1-r<\Vert z\Vert <1\}$, $\Vert \cdot
 \Vert $ is a polydisk norm in $\cc^q$.
 The following lemma for $q=1$ can be found in [I$_4$], Lemma 2.2.1.
 Proof for any $q\ge 1$ is the same.

 \smallskip\noindent\bf
 Lemma 1. \it If any meromorphic map $f:E_n^q(r)\to Y$
 extends to a meromorphic map $\hat f:\Delta^{n+q}\to Y $ then the space
 $Y$ possesses a meromorphic extension property in bidimension $(n,q)$.
 
 \smallskip\rm
 We shall make use also from one result on \it meromorphic families of
 analytic subsets \rm from [I$_4$].

 \medskip \rm
 Let  $S$ be a set, and $W\subset\subset \cc^q$ an open subset. $W$ is equipped
 with the usual Euklidean metric from $\cc^q$. $Y$ is again some complex space.

 \smallskip \noindent
 \bf Definition. \it (i) By a family of $q$-dimensional analytic subsets in
 complex
 space $X=W\times Y$ we shall understand an subset ${\cal F}\subset S\times 
 W\times X $ such that, for every $s\in S$ the set  $ {\cal F}_s = {\cal F}
 \cap \{s\} \times W\times X$ is a graph of a meromorphic mapping of $W$
 into $X$.

 \smallskip

 (ii) If the  set $S$ is equipped with topology and the space $X$ is equipped
 with some Hermitian metric $h$ we say that the family ${\cal F}$ is
 continuous at  point $s_0\in S$ if ${\cal H}-\lim_{s\rightarrow s_0}
 {\cal F}_s = {\cal F}_{s_0}$.

 \smallskip
 (iii) When $S$ is a complex space itself call the family ${\cal F}$
 meromorphic if the closure $\hat {\cal F}$ of the set ${\cal F}$ is an
 analytic subset of $S\times W \times X$.

\smallskip\rm
 Here by ${\cal H}-\lim_{s\rightarrow s_0} {\cal F}_s $
 we denote the limit of closed subsets of ${\cal F}_s$ in the Hausdorff
 metric on $W\times X$. ${\cal F}$ is continuous if it is
 continuous
 at each point of $S$. If $W_0 $ is open in $W$ then the
 restriction ${\cal F}_{W_0}$ is naturally defined as ${\cal F}\cap (S
 \times W_0 \times X)$.

The statement about meromorphic families we need can be formulated as
follows. For the standart notions and facts from pluripotential theory
we refer to [Kl].

Consider a
meromorphic mapping $f: V \times W_0 \longrightarrow X$ into a
complex space $X$, where $V$ is a domain in $\cc^p$. Let $S$ be some closed
subset of $V $ and $s_0\in S$
some accumulation point of $S$. Suppose that for each $s\in S$ the
restriction $f_s = f|_{\{ s\} \times W_0} $ meromorphically extends
onto $W\supset\supset W_0$. We suppose additionally that there is a
compact $K\subset\subset X$ such that for all $s\in S$ $f_s(W)\subset K$.

Let $\nu _j$ denotes the minima of
volumes of $j$-dimensional compact analytic subsets contained in our compact
$K\subset X$. $\nu_j>0$, see Lemma 2.3.1 from [I$_4$]. Fix some
$W_0\subset\subset W_1\subset\subset W$ and put

\smallskip

$$
\nu = \min \{ {\sl vol}(A_{q-j})\cdot \nu _j : j = 1,\ldots,q\},\eqno(3)
$$

\smallskip\noindent
where $A_{q-j}$ are running over all $(q-j)$-dimensional analytic subsets of
$W$, intersecting $\bar W_1$.  Clearly
$\nu > 0$. In the
following {\sl Lemma } the volumes of graphs over $W$
are taken. More precisely, having an Euklidean metric form $w_e=dd^c\Vert
z\Vert^2$ on $W\subset \cc^q$ and Hermitian metric form $w_h$ on $X$, we
consider $\Gamma_{f_s}$ for $s\in S$ as an analytic subsets of $W\times X$
and their volumes are

$$
vol(\Gamma_{f_s}) = \int_{\Gamma_{f_s}}(p_1^*w_e+p_2^*w_h)^q =
\int_W(w_e+(p_1)_*p_2^*w_h)^q,\eqno(4)
$$

\noindent
where $p_1:W\times X\to W$ and $p_2:W\times X\to X$ are natural projections.

\smallskip\noindent
\bf Lemma 2. \it Suppose  that there exists
a neighbourhood $U\ni s_0$ in $V$ such that, for all $s_1,s_2\in S\cap U$

\smallskip

$$
| {\sl vol}(\Gamma _{f_{s_1}}) - {\sl vol}(\Gamma _{f_{s_2}})
| < \nu /2. \eqno(5)
$$

\smallskip\noindent
If  $s_0$ is a locally regular point of $S$ then  there exists
a neighbourhood $V_1\ni s_0$ in $V$ such, that $f$
meromorphically extends onto $V_1\times W_1$.

\smallskip\rm
 Further, slightly modifying arguments from [I$_3$] we shall derive now the
 following version of so called Continuity principle.

 Let   $f:H_n^q(r)\longrightarrow Y$ be a
 given  meromorphic mapping. Let
 $A_s^q(1-r,1):= \{ s\}\times A^q(1-r,1)$ for $s\in \Delta^n$. We suppose
 that for $s$ in some nonempty subset $S\subset \Delta^n$ the restriction
 $f_s:=f\mid_{A^q_s(1-r,1)}$ is well defined and extends meromorphically to
 the polydisc $\Delta^q$. .

 \smallskip\noindent\bf
 Lemma 3. \it Suppose that $f:H_1^q(r)\longrightarrow Y$ is
 meromorphic and:

 (i) there is a  compact $K\subset\subset Y$ such that $f(\Delta^1
 \times A^q(1-r,1))\subset K$ and $f(\Delta_s^q)\subset K$ for all $s\in S$;

 (ii) there is a
 constant $C_<0 \infty $ such that ${\sl vol}(\Gamma_{f_s})\le C_0$  for all
 $s\in S$.

 \noindent Then:

 1. Either there is a neighborhood $U\ni 0$ in $\Delta^n$
 and a meromorphic extension of $f$ onto  $U\times \Delta^q$, or

 2. $0$ is an isolated point of $S$.

 \smallskip\rm The volumes here are measured with respect to the Euclidean
 metric on $\bf C^q$ \rm and some Hermitian metric $h$ on $Y$. The condition
 of boundedness in (ii) clearly  does not depend on the particular choice of
 $h$. We shall refer to this statement as to C.P. The condition $n=1$ is
 important here, see Example 1 in [I$_3$]. We shall also discuss the related
 questions in our problem section.

 \smallskip\rm
 To derive the proof of this statement from the reasonings in [I$_3$] we shall
 need some notions and results from the theory of cycle spaces
 (due to  D.Barlet, see [B$_2$]) as they where adapted to our "noncompact"
 situation in [I$_3$]. For the english spelling of the Barlet terminology we
send an interested reader to [Fj].

Recall that an analytic cycle of dimension $q$ in complex space $Y$ is a
formal sum $Z=
\sum_jn_jZ_j$, where $\{ Z_j\} $ is a locally finite sequence of  analytic
subsets (allwayse of pure dimension
$q$) and $n_j$ are positive integers called multiplicities of $Z_j$.
$\vert Z\vert :=\bigcup_jZ_j$-support of $Z$. All complex spaces in this
paper are reduced, normal and countable at infinity.

\smallskip
With a given  meromorphic mapping $f:\Delta\times A^q(1-r,1)\to X$, satifying
conditions of Lemma 3 we associate the following space
of cycles. Fix some $0<c<1$.

\smallskip
Consider a set ${\cal C}^{'}_{f,C}$ of all
analytic cycles $Z$ in $Y:=\Delta
^{1+q}\times X$ of pure dimension $q$, such that:

(a) $Z\cap [\Delta\times A^q(1-r,1]\cap X] = \Gamma_{f_{z}}\cap \{ z\}
\times A^q(1-r,1]\times X$ for
some $z\in \Delta(c)$. This means, in particular, that for this $z$
mapping $f_{z}$ extends meromorphically from $A^q_{z}(1-r,1)$ onto
$\Delta^q_{z}$.

(b) ${\sl vol}(Z)< C$, where $C$ is a some constant, $C>C_0$, $C_0$
beeing from Lemma 3.

\smallskip
Define $\bar {\cal C}_{f,C}$ to be a closure of ${\cal C}^{'}_{f,C}$
in the usual topology of currents, see below. In [I$_3$] it was
 shawn that ${\cal C}
_{f,C}:=\{ Z\in \bar {\cal C}_{f,C}: {\sl vol}(Z)<C\} $
is an analytic space of finite dimension in the neighborhood of each
of its points.

As we already had mentioned our first aime is to prove the analyticity of
${\cal C}_{f,C} $.

Let $f:\Delta\times A^q(1-r,1]\to X$ be our map. Denote by ${\cal C}_0$ the
subset of $\bar {\cal C}_{f,C}$ consisting of cycles which are limits of
$\{ \Gamma_{f_{s_n}}\} $ for $s_n\rightarrow 0,s_n\in S$. This is a compact 
subset (by Bishop's theorem) of the topological space ${\cal C}_{f,2C}$. For 
every cycle $Z\in {\cal C}_0$ define its neighborhood $W_Z$ as above. Let 
$W_{Z_1},...,W_{Z_N}$ be
a finite covering of ${\cal C}_0$. Remark that there is an $\eps_0>0$ such that
for any $s\in S\cap \Delta(\eps_0)$ we have 
$\Gamma_{f_s}\subset \bigcup_{j=1}^NW_{Z_j}$.

 Now we are prepared to sketch the proof of Lemma 3.
 Consider a universal family ${\cal Z}:=\{ Z_a:a\in {\cal C}_{f,2C_0}\} $.
 This is
 complex space of finite dimension . We have an evaluation map

 $$
 F:{\cal Z}\to \Delta^{1+q}\times X
 $$

 \noindent
 defined by $Z_a\in {\cal Z}\to Z_a\subset \Delta^{1+q}\times X$.
 Consider the union $\hat {\cal C}_0$ of those components
 of ${\cal C}_{f,2C_0}$ which intersect ${\cal C}_0$. Recall, that ${\cal
 C}_0$ stands here for the set of all limits of $\{ \Gamma_{f_{s_n}}, s_n\in
 S\}$. At least one of those
 components, say ${\cal K}$, containes two points $s_1$ and $s_2$ s.t.
 $Z_{s_1}$ projects onto $\Delta^k_0$ and $Z_{s_2}$ projects onto
 $\Delta^k_s$ with $s\not= 0$. This is just because $S$ contains more then
 one point. Consider the restriction ${\cal Z}\mid_{{\cal K}}$ of the
 universal space onto ${\cal K}$. This is an irreducible
 complex space of  finite dimension. Take points $z_1\in Z_{s_1}$
 and $z_2\in Z_{s_2}$ and join them by an analytic disk $\phi :\Delta \to
 {\cal Z}\mid_{{\cal K}}$, $\phi (0)=z_1, \phi (1/2)=z_2$. Then the
 composition $\psi = \pi\circ F\circ \phi :\Delta \to \Delta $ is not
 degenerate because $\psi (0)=0\not= s=\psi (1/2)$. Thus $\psi $ is proper
 and obviously so is the map $F:{\cal Z}\mid_{\phi
 (\Delta )}\to F({\cal Z}\mid_{\phi (\Delta )})\subset \Delta^{1+q}\times
 X$. Thus  $F({\cal Z}\mid_{\phi (\Delta )})$ is an analytic set in $U\times
 \Delta^k\times X$ for small enough $U$ extending $\Gamma_f$ by the reason
 of dimension.

 \smallskip
 We also shall make use
 also of the following result due to D. Barlet, see [B$_1$] Proposition 3:

 \medskip\noindent\bf
 Lemma 4. \it Let $X$ be a reduced complex space (of finite
 dimension) and let $\rho :X\longrightarrow \bf R^+ $ \it be a strictly 
 $q$-convex
  function. Let $h$ be some $C^2$-smooth Hermitian metric on $X$. Then there
 exists an Hermitian metric $h_1$ and a function $c:\bf R^+\longrightarrow
 R^+$ \it (both of class $C^2$) such that:

 (i) $h_1 \ge h$;

 (ii) the $(q,q)$ -  form $\Omega=dd^c[(c\circ\rho)w_{h_1}^{q-1}]$ is
 strictly positive  on $X$.

 \smallskip\rm Here $w_h$ is the (1,1)-form canonically associated with $h$.
 In our case we need $X=\Delta^{n+q}\times Y$ and we shall use only  the
 fact that on $X$ there exists a strictly positive $(q,q)$-form which is
 $dd^c$-exact:  in fact $d$-exactness is going to be sufficient for us.
 We recall that a $(q,q)$-form $\Omega $ is called strictly positive if for
 any $x\in X$ and linearly independent vectors $v_1,...,v_q\in T_xX$ one has
 $\Omega_x(iv_1\wedge\bar v_1,...,iv_q\wedge\bar v_q)>0$.

 \bigskip\noindent\bf 2. Proof of the Theorem. \rm
 
 \smallskip\noindent\sl
 Step 1. Case $n=1$. \rm

 \smallskip

 Let  $f:H_1^q(r) \longrightarrow Y$ be our
 meromorphic mapping. Let us denote by $W$ the biggest open subset of
 $\Delta^1$ such that $f$ extends meromorphically  to $H_W(r):=
 (\Delta^1\times A^q(1-r,1)) \cup (W\times\Delta^q)$, and let us remark
 explicitly that the  complex space  $X=\Delta^{1+q}\times Y$ is (obviously)
 $q$-complete.

 We apply Barlet's
 Theorem, see $\S 1$, by taking as $\rho$ a strictly $q$-convex exhaustion of
 $X$ in order to have a strictly positive $dd^c$-exact  $(q,q)$-form
 $\Omega $ on $X$. Let $w$ be a fixed  $(q-1,q-1)$-form of class $C^2$ such
 that $dd^cw=\Omega$. Let us denote by $F$ the set of points of
 indeterminancy of $f$.

 By shrinking the polydisc $\Delta^{1+q}$, we can suppose, withought loss of
 generality, that $f_z$ is defined in the neighborhood of $\bar\Delta^q$ for
 all $z\in W $. In the same way we can suppose that $w\in C^2(\bar\Delta^{1+q}
 \times Y)$, i.e. is smooth up to the boundary.

 We need to prove that $W =\Delta $. Suppose not, and fix a point $z_0\in
 \partial W \cap  \Delta $. Denote by $V$ some disc centered at
 $z_0$ which whichis contained in $\Delta $. For $z\in V\cap W $
  one has

  $$
  {\sl vol}(\Gamma_{f_z})=\int_{\Gamma_{f_z}}\Omega=\int_{\Gamma_{f\mid_
  {\partial\Delta^q_z}}}d^cw\le C\eqno(6)
  $$

  \noindent
  where the constant $C$ does not depend on $z\in V\cap W $, while $d^cw$ iso
  f class $C^1$ on $\bar\Delta^{1+q}\times Y$.

  To obtain estimate $(1)$ we
  had used the  fact that we can  measure the volumes of analytic sets of
  pure dimension $q$ contained in some compact part  of $X$ by means of
  $\int \Omega $ with $\Omega$ a strictly positive $(q,q)$-form   on
  $X$.

  We are going to check if the conditions  of the Continuity Principle,
  mentioned in $\S 1$ are satisfied. The inequality (6) says that the second
  assumption of C.P. is satisfied.
  To  check if the first one is satisfied,
   let  us suppose that there exists a sequence $\{ z_n\} \subset  V
   \cap \Omega $,
  converging to $z_{\infty }\in \Delta $, such that $\{ \Gamma_{\nu }:=
  \Gamma_{f_{z_{\nu }}}\} $ is not contained in any relatively compact subset
  of $\bar\Delta^{1+q}\times Y$. If  $\nu $ is big enough, the restriction
  $\rho\mid_{\Gamma_{\nu }} $ will have then a strict maximum  in the
  interior of $\Gamma_{\nu } $. This is impossible because the Levi form of
  $\rho\mid_{\Gamma_{\nu }}$ has at least one positive eigenvalue at each
  point of $\Gamma_{\nu }$. Remark also that $q$-complete space $Y$ cannot
  containe any compact $q$-dimensional subspace.
  C.P.  says now (since $W$
  is not contained in any proper analytic subset of  any  neighborhood of
  $z_0$)  that $f$ meromorphically extends to $V_1\times \Delta^q$ for some
  neighborhood $V_1$ of $z_0$ in $\Delta $. This proves that $W = \Delta $.

  \smallskip\noindent\sl
  Step 2. Case $n\ge 2$. \rm

  \smallskip
  This will done by induction on $n$. As is follows from Lemma 1, all
  we need is to extend the mappings from $E_n^q(r)$ to $\Delta^{n+q}$.
  For $n=1, E_1^q(r)=H_1^q(r)$ and thus this is already done by  
  Step 1.

  Notice that $E_{n+1}^q(r)=\Delta \times E_n^q(r)$. Denote
  $E_{n,z}^q(r):=\{ z\} \times E_n^q(r)$ for $z\in \Delta $. Remark
  that by the induction  hyposethis the restriciton $f\mid_{E_{n,z}^q(r)}$
  meromorphically extends onto $\Delta^{n+q}_z:=\{ z\} \times \Delta^{n+q}$
  for all $z\in \Delta $. We denote by $W $ the maximal open subset in
  $\Delta $ such  that our map $f$ extends meromorphically onto $W\times
  \Delta^{n+q}$.

  Put $S=\Delta \setminus W $ and consider a family $\{ \Gamma_{f_s}:
  s\in S\}
  $ of analytic subsets in $X:=\Delta^{n+q}\times Y$. Here, as usually
  by $\Gamma_{f_s}$ we denote the graph of the restriction $f_s:=f\mid_{
  \Delta^{n+q}_s}$.
  Define $S_k:= \{ s\in S: vol(\Gamma_{f_s})\le k\cdot {\nu \over 2}\} $.
  Where $\nu $ is from Lemma 2 with $W=\Delta^{n+q}, W_0=\Delta^{n+q}_{1-r/2}$.
  By maximality of $S$ and by Lemma 2 we see that all points of each $S_k$
  are locally regular, thus each $S_k$ is polar. So $S$ is a polar subset of
  $\Delta $, in other words it is a set of harmonic measure zero in $\Delta $.

  By some linear coordinate transformation in
  $\bf C^{1+n+q}$ we are going to change  a little bit the band of the
  $\Delta^{n+q}$-direction, in order to prove in the same manner that $f$
  meromorphically extends to the whole of $\Delta^{1+n+q}$.   In fact, let
  us consider linear changes $L$ of the coordinate system in $\bf C^{1+n+q}$
  \rm  whose associated matrices are of the form $(L_1,L_2)$, where $L_1$ is
  (a number) close to zero and $L_2$ is close to the identity map of
  $\cc^{n+q}$ into itself. For each $L$ of this   form we can extend $f$
  onto $\Delta^{1+n+q}\setminus L^{-1}(S^L\times \Delta^{n+q})$,             
  where $S^L$ is a  set of harmonic measure zero in $\Delta $.                                                         
  $ \Sigma $ in appropriate coordinate system is a $(1+n+q)$-product of
  the closed sets of harmonic measure zero on the plain. Thus $\Sigma $
  is pluripolar and of Haussdorff dimension zero.
  Using the fact that $\Delta^{1+n+q}\times Y$ is (obviously!) $(n+q)$ -
  complete and Lemma 3 we can remove the singularity $\Sigma $.

  \bigskip

  \hfill{q.e.d.}

  \bigskip\noindent\bf 3. Consequences and open questions.
  
  Let us start with some direct consequences of the Theorem.

  \medskip\noindent\bf
  Corollary 1. \rm (Thullen type extension Theorem) \it Let $\Omega\subset\bf
   C^n$ \it be an open subset, $V\subset \Omega$
  be an analytic subvariety of dimension $q$ and $G$ be an open subset of
  $\Omega$ which intersects every $q$-dimensional branch of $V.$ Every
  meromorphic mapping $f:(\Omega\setminus V)\cup G\longrightarrow Y$,
  where $Y$ is a $q$-complete complex space, extends to a meromorphic
  mapping from $\Omega$ to $Y.$\rm

  \smallskip\rm
  In particular one has
  
  \medskip\noindent\bf
  Corollary 2. \rm (Riemann type extension Theorem) \it Let
  $\Omega\subset\bf C^n$ \it be an open subset, $V\subset \Omega$ be an
  analytic subvariety of dimension $q-1$. Every meromorphic mapping
  $f:\Omega\setminus V\longrightarrow Y$, where $Y$ is a $q$-complete
  complex space, extends to a meromorphic mapping from $\Omega$ to $Y.$
  
  \rm The proofs of Corollaries 1 and 2 are immediate after [S$_1$],
  p.5.

  A general Thullen type extension Theorem for meromorphic mappings is
  proved by Siu when $Y$ is compact K\"ahler in [S$_2$].  We have also:

  \medskip\noindent\bf
  Corollary 3. \it Let $Y$ be a   complex analytic  space of dimension $q$
  and let us suppose that every irreducible component
  of $Y$ of dimension $q$ is non compact.  Every meromorphic mapping
  $f:H_n^q(r)\longrightarrow Y$, extends to a meromorphic mapping from
  $\Delta^{n+q}$ to $Y$.

  \rm\medskip In fact, every complex space of dimension $n$ with no
  compact irreducible component of dimension $n$ is $n$-complete, by a
  Theorem of Ohsawa, [O] Th. 1.

  \bigskip
  We shall end with discussing some open questions, which naturally arise
  from the results and attempts of this paper.

  \smallskip\noindent\sl
  Question 1. \it Let $Y$ be a compact complex three-fold. Prove that
  every meromorphic (or holomorphic) map $f:H_1^2(r)\to Y$ extends onto
  $\Delta^3\setminus \{ discrete$ $set$ $of$ $points\} $.

  In particular, if $K\subset\subset \Delta^3$ with connected complement
  then every meromorphic map $f\Delta^3\setminus K\to Y$ extends onto
  $\Delta^3$ minus finite set of points.

  \smallskip\rm
   For the proof of such type of statements one can try to use special
metrics on $Y$. Namely a compact complex three-fold possesses a Hermitian
metric $h$, such that its associated $(1,1)$-form $\omega_h$ satisfies
$dd^c\omega_h^2=0$. This can help to bound the volumes of the images of
two-disks in $Y$.

  \medskip
  The next question comes out when one tries to prove the Corollary
  from the Introduction withough assuming $Y$ to be noncompact.

  \smallskip\noindent\sl
  Question 2. \it Let $Y$ is a compact complex manifold (space) of dimension
  $q\ge 2$. Suppose that there exists a meromorphic map $f:\bf B_*^{q+1}\rm
  \to Y$ from punctured ball in $\cc^{q+1}$ onto $Y$ such that for any
  $\eps >0$ the restriction $f_{\eps }:= f\mid_{\bf B^{q+1}_*(\eps )}$
  of $f$ onto the punctured $\eps $-ball is still surjective. Prove that
  $Y$ is Moishezon.

  \smallskip\rm  In the case of positive answer to this question, one can
extend this $f$ meromorphically to zero.

  In fact one need somewhat stronger stetment. Let $M\ni 0$ be a strongly
pseudoconvex hypersurface in the ball ${\bf B_*^{q+1}}$, which divides it
into two parts $B^+$ and $B^-$. Let a meromorphic map $f:B^+\to Y$ as in
Question 2 is given. Suppose that $M$ is concave from the side of $B^+$.

\smallskip\it Prove that if for any $\eps >0$ the restriction
$f\mid_{\bf B_{\eps }^{q+1}\cap B^+}:\bf B_{\eps }^{q+1}\cap B^+\to Y$ \it is
surjective, then $f$ extends meromorphically to zero.

  \smallskip\noindent\sl
  Question 3. \it Can one remove the condition $n=1$ from Lemma 3?

  \smallskip\rm
   The attepts leed to the "non analytic" version of Remmert
   proper mapping theorem and to the questions of locall flattenings.
   The major problem here is that ${\cal C}_{f,C}$ will be not
an analytic space now and $F$ will be not proper in general.

  \bigskip

  \centerline{ \bf R E F E R E N C E S}

  \bigskip\noindent\rm [B$_1$] Barlet D.:\it Convexit\'e de l'espace des

  cycles.
  \rm Bull. Soc. Math. France, \bf  106 \rm (1978) 373-397

  \smallskip\noindent [B$_2$] Barlet D.: \it Espace analytique reduit de
  cycles analytiques \- complexers compacts d'un espace analytique complexe
  de dimension finie. Seminair Norguet
  IX, Lect. Notes Math., 1-157 (1975).

  \smallskip\noindent [Fj] Fujuki A.: Closedness of the Douady Space of
Compact K\"ahler Spaces. Publ. RIMS, Kyoto Univ. {\bf 14}, 1-52 (1978).

  \smallskip\noindent [H] Hartogs F.:\it Zur Theorie der analytischen
  Funktionen me\-hrerer un\-abh\"angiger
  
  \noindent Ver\"andenlichen insbesondere
  \"uber die Darstellung derselben durch Reihen, wel\- che nach Potenzen
  einer Ver\"anderlichen fortschreiten.\rm  Math. Ann.\bf 62  \rm (1906) 1-88

  \smallskip\noindent [I$_1$] Ivashkovich S.: \it The Hartogs type extension
  Theorem for meromorphic maps into compact K\"ahler manifolds. \rm Invent.
  Math.\bf 109 \rm (1992) 47-54.

  \smallskip\noindent [I$_2$] Ivashkovich S.: \it Spherical shells as
  obstructions for the extension of holomorphic mappings. \rm The Journal
  of Geometric Analysis.\bf 2 \rm (1992) 683-692

  \smallskip\noindent [I$_3$] Ivashkovich S.:\it Continuity principle
  and extension properties of meromorphic mappings with values in non
  K\"ahler manifolds. \rm MSRI Preprint No. 1997-033. math.CV/9704219

  \smallskip\noindent [I$_4$] Ivashkovich S.:\it One example in concern
  with extension and separate analyticity properties of meromorphic
  mappings. \rm To appear in Amer. J. Math. math.CV/9804009

  \smallskip\noindent [Kl] Klimek M.: \it Pluripotential theory. \rm
  London. Math. Soc. Monographs, New Series 6 (1991).

  \smallskip\noindent [O] Ohsawa T.: \it Completeness of Noncompact Analytic
  Spaces. \rm Publ. RIMS Kyoto Univ.\bf 20 \rm (1984) 683-692

  \smallskip\noindent [R] Remmert R.:\it Holomorphe und meromorphe
  Abbildungen komplexer R\"aume.
  \rm Math. Ann.\bf 133 \rm (1957) 328-370

  \smallskip\noindent [S$_1$] Siu Y.-T.: \it Techniques of extension of
  analytic objects.\rm Marcel Dekker, New York , (1974).

  \smallskip\noindent [S$_2$] Siu Y.-T.: \it Extension of meromorphic maps
  into  K\"ahler manifolds. \rm Ann.Math. \bf 102  \rm (1975) 421-462

  \smallskip\noindent [St] Stein K.: Rocky Montain Math. J. (1972)

  \bigskip\noindent Sergei Ivashkovich - Universite' des Sciences et
  Technologies de Lille - UFR de Mathematiques - 59655 Villeneuve d'Ascq
  Cedex (France) - email: ivachkov@gat.univ-lille1.fr

  \bigskip
  \noindent Alessandro Silva - Universita' di Roma La Sapienza - Dipartimento
   di Matematica
  G. Castelnuovo - P.le A.Moro - 00185 Roma (Italy) - email:
   silva@mat.uniroma1.it

   \end